\documentclass[11pt]{amsart}

\usepackage{amssymb}
\usepackage[colorlinks]{hyperref}
\usepackage[margin=1.24in]{geometry}
\usepackage{enumitem}
\usepackage{verbatim}
\usepackage{tikz}
\usetikzlibrary{matrix,arrows,decorations.markings}
\usepackage{caption}

\theoremstyle{definition}

\DeclareMathOperator{\Spec}{Spec}
\DeclareMathOperator{\Hom}{Hom}

\title{Computing with jets}
\date{\today}

\author{Federico Galetto}
\address{Department of Mathematics and Statistics, Cleveland State University, Cleveland, OH, 44115-2215, USA}
\email{\href{mailto:f.galetto@csuohio.edu}{\nolinkurl{f.galetto@csuohio.edu}}}
\urladdr{\href{https://math.galetto.org}{\nolinkurl{https://math.galetto.org}}}

\author{Nicholas P. Iammarino}
\address{Department of Mathematics and Statistics, Cleveland State University, Cleveland, OH, 44115-2215, USA}
\email{\href{mailto:nickiammarino@gmail.com}{\nolinkurl{nickiammarino@gmail.com}}}



\begin{document}

\begin{abstract}
  We introduce a Macaulay2 package for working with jet schemes. The
  main method constructs jets of ideals, polynomial rings and their
  quotients, ring homomorphisms, affine varieties, and
  (hyper)graphs. The package also includes additional methods to
  compute principal components and radicals of jets of monomial
  ideals.
\end{abstract}

\maketitle{}

\section{Introduction}
\label{sec:introduction}

Roughly speaking, the scheme of $s$-jets of a scheme $X$ is the
collection of order $s$ Taylor approximations at points of $X$.  More
formally, let $X$ be a scheme over a field $\Bbbk$. Following \cite[\S
2]{MR2483946}, we call a scheme $\mathcal{J}_s (X)$ over $\Bbbk$ the
scheme of $s$-jets of $X$, if for every $\Bbbk$-algebra $A$ there is a
functorial bijection
\begin{equation*}
  \Hom ( \Spec (A), \mathcal{J}_s (X) ) \cong
  \Hom ( \Spec (A[t]/\langle t^{s+1} \rangle), X ).
\end{equation*}
This means that the $A$-points of $\mathcal{J}_s (X)$ are in bijection
with the $A[t]/\langle t^{s+1} \rangle$-points of $X$. It follows that
$\mathcal{J}_0 (X) \cong X$, and $\mathcal{J}_1 (X)$ is the total
tangent scheme of $X$, in line with the definition of tangent space
using dual numbers \cite[II, Exercise 2.8]{MR0463157}. Jet schemes
play an important role in the study of singularities, as initially
suggested by J.~Nash \cite{MR1381967}, and in connection with other
related topics, such as motivic integration and birational geometry
\cite{MR1905328,MR1856396,MR1896234,MR2483946}.

The existence of jet schemes is proved in detail in \cite[\S
2]{MR2483946}. We recall an essential step, which is the construction
of jets of an affine variety. Let $X$ be an affine variety over
$\Bbbk$. Consider a closed embedding of $X$ into an affine space
$\mathbb{A}^n$ over $\Bbbk$. Let $I = \langle f_1,\dots, f_r\rangle$
be the ideal of $R = \Bbbk [x_1,\dots,x_n]$ corresponding to this
embedding. For $s\in \mathbb{N}$, define the polynomial ring
\begin{equation*}
  \mathcal{J}_s (R) = \Bbbk [x_{i,j} \,|\, i=1,\dots,n, j=0,\dots,s].
\end{equation*}
For each $k=1,\dots,n$, perform the substitution
\begin{equation*}
  x_k \mapsto x_{k,0} + x_{k,1} t + x_{k,2} t^2 + \dots + x_{k,s} t^s
  = \sum_{j=0}^s x_{k,j} t^j
\end{equation*}
taking elements of $R$ to elements of $\mathcal{J}_s (R)
[t]$. Applying this substitution to a generator $f_i$ of $I$ gives the
following decomposition:
\begin{equation*}
  f_i \left( \sum_{j=0}^s x_{1,j} t^j, \dots, \sum_{j=0}^s x_{n,j} t^j \right)
  = \sum_{j\geqslant 0} f_{i,j} t^j,
\end{equation*}
where the coefficients $f_{i,j}$ are polynomials in
$\mathcal{J}_s (R)$. The \emph{ideal of $s$-jets} of
$I = \langle f_1,\dots,f_r\rangle$ is the ideal of $\mathcal{J}_s (R)$
defined by
\begin{equation*}
  \mathcal{J}_s (I) = \langle f_{i,j} \,|\,
  i=1,\dots,r, j=0,\dots,s\rangle.
\end{equation*}
The scheme of $s$-jets of $X$ is
$\Spec (\mathcal{J}_s (R) / \mathcal{J}_s (I) )$.

This paper introduces the \texttt{Jets} package\footnote{Available at
  \url{https://github.com/galettof/Jets}.} for Macaulay2 \cite{M2},
streamlining the process of constructing ideals of jets as indicated
above. We adopt the following notation: the variables in the
polynomial rings containing the equations of jets have the names of
the variables of the original equations with the order of the jets
appended to them, and the same subscripts. Moreover, the rings
containing the equations of jets are constructed incrementally as
towers.

Ideals of jets are computed via the \texttt{jets} method applied to
objects of type \texttt{Ideal}. In addition, the \texttt{jets} method
can also be applied to objects of type \texttt{QuotientRing},
\texttt{RingMap}, and \texttt{AffineVariety}, with the effects one
would expect from applying jet functors.  For more information,
including grading options, we invite the reader to consult the
documentation of the package. We showcase some package applications
below.

We are grateful to Greg Smith for valuable feedback on an early
version of the package.

\section{Jets of monomial ideals}
\label{sec:jets-monomial-ideals}

As observed in \cite{MR2229478}, the ideal of jets of a monomial ideal
is typically not a monomial ideal.

\begin{verbatim}
i1 : needsPackage "Jets";

i2 : R=QQ[x,y,z];

i3 : I=ideal(x*y*z);

o3 : Ideal of R

i4 : J2I=jets(2,I);

o4 : Ideal of QQ[x0, y0, z0][x1, y1, z1][x2, y2, z2]

i5 : netList J2I_*

     +---------------------------------------------------------------+
o5 = |y0*z0*x2 + x0*z0*y2 + x0*y0*z2 + z0*x1*y1 + y0*x1*z1 + x0*y1*z1|
     +---------------------------------------------------------------+
     |y0*z0*x1 + x0*z0*y1 + x0*y0*z1                                 |
     +---------------------------------------------------------------+
     |x0*y0*z0                                                       |
     +---------------------------------------------------------------+
\end{verbatim}

However, by \cite[Theorem 3.1]{MR2229478}, the radical is always a
(squarefree) monomial ideal. In fact, the proof of \cite[Theorem
3.2]{MR2229478} shows that the radical is generated by the terms of
the jet equations constructed as in the introduction. This observation
provides an alternative algorithm for computing radicals of jets of
monomial ideals, which can be faster than the default radical
computation in Macaulay2.

\begin{verbatim}
i6 : jetsRadical(2,I);

o6 : Ideal of QQ[x0, y0, z0][x1, y1, z1][x2, y2, z2]

i7 : netList pack(5,oo_*)

     +--------+--------+--------+--------+--------+
o7 = |y0*z0*x2|x0*z0*y2|x0*y0*z2|z0*x1*y1|y0*x1*z1|
     +--------+--------+--------+--------+--------+
     |x0*y1*z1|y0*z0*x1|x0*z0*y1|x0*y0*z1|x0*y0*z0|
     +--------+--------+--------+--------+--------+
\end{verbatim}

For a monomial hypersurface, \cite[Theorem 3.2]{MR2229478} describes
the minimal primes of the ideal of jets. Moreover, the main theorem in
\cite{math/0607638} counts the multiplicity of the jet scheme of a
monomial hypersurface along its minimal primes (see also
\cite{MR2372306}). We compute the minimal primes, then use the
\texttt{LocalRings} package \cite{LocalRingsSource} to compute their
multiplicities in the second jet scheme of the example above.

\begin{verbatim}
i8 : P=minimalPrimes J2I;

i9 : --flatten ring to use LocalRings package
     (A,f)=flattenRing ring J2I;

i10 : needsPackage "LocalRings";

i11 : --quotient by jets ideal as a module
      M=cokernel gens f J2I;

i12 : --compute the multiplicity of the jets along each component
      mult=for p in P list (
          Rp := localRing(A,f p);
          length(M ** Rp)
          );

i13 : netList(pack(4,mingle{P,mult}),HorizontalSpace=>1)

      +--------------------+---+--------------------+---+
o13 = | ideal (z0, y0, x0) | 6 | ideal (z0, y0, z1) | 3 |
      +--------------------+---+--------------------+---+
      | ideal (z0, y0, y1) | 3 | ideal (z0, x0, z1) | 3 |
      +--------------------+---+--------------------+---+
      | ideal (z0, x0, x1) | 3 | ideal (z0, z1, z2) | 1 |
      +--------------------+---+--------------------+---+
      | ideal (y0, x0, y1) | 3 | ideal (y0, x0, x1) | 3 |
      +--------------------+---+--------------------+---+
      | ideal (y0, y1, y2) | 1 | ideal (x0, x1, x2) | 1 |
      +--------------------+---+--------------------+---+
\end{verbatim}

\section{Jets of graphs}
\label{sec:jets-graphs}

Jets of graphs were introduced in \cite{2104.08933}. Starting with a
finite, simple graph $G$, one may construct a quadratic squarefree
monomial ideal $I(G)$ (known as the \emph{edge ideal} of the graph) by
converting edges to monomials (see for example \cite{MR3184120}).  One
may then consider the radical of the ideal of $s$-jets of $I(G)$,
which is again a quadratic squarefree monomial ideal. The graph
corresponding to this ideal is the graph of $s$-jets of $G$, denoted
$\mathcal{J}_s (G)$.

Jets of graphs and hypergraphs can be obtained by applying the
\texttt{jets} method to objects of type \texttt{Graph} and
\texttt{HyperGraph} from the Macaulay2 \texttt{EdgeIdeals} package
\cite{EdgeIdealsSource,MR2878668} (which is automatically loaded by
the \texttt{Jets} package). Consider, for example, the graph in Figure
\ref{fig:1}.

\begin{figure}[htb]
  \centering
  \begin{tikzpicture}[scale=0.5]
    % styles
    \tikzset{vertex/.style = {shape=circle,draw,thick,inner sep=2pt}}
    \tikzset{edge/.style = {thick}}
    \tikzset{every label/.append style={font=\small}}

    % vertices
    \node[vertex] (a) at (90:2) [label=90:{$a$}] {};
    \node[vertex] (b) at (18:2) [label=18:{$b$}] {};
    \node[vertex] (c) at (306:2) [label=306:{$c$}] {};
    \node[vertex] (d) at (234:2) [label=234:{$d$}] {};
    \node[vertex] (e) at (162:2) [label=162:{$e$}] {};

    % edges
    \draw[edge] (a) -- (c) -- (b) -- (e) -- (c);
    \draw[edge] (e) -- (a) -- (d) -- (b);
  \end{tikzpicture}

  \caption{}\label{fig:1}
\end{figure}

\begin{verbatim}
i1 : needsPackage "Jets";

i2 : R=QQ[a..e];

i3 : G=graph({{a,c},{a,d},{a,e},{b,c},{b,d},{b,e},{c,e}});
\end{verbatim}

We compute the first and second order jets, and list their edges.

\begin{verbatim}
i4 : J1G=jets(1,G); netList pack(7,edges J1G)

     +--------+--------+--------+--------+--------+--------+--------+
o5 = |{c1, a0}|{d1, a0}|{e1, a0}|{c1, b0}|{d1, b0}|{e1, b0}|{a1, c0}|
     +--------+--------+--------+--------+--------+--------+--------+
     |{b1, c0}|{e1, c0}|{a0, c0}|{b0, c0}|{a1, d0}|{b1, d0}|{a0, d0}|
     +--------+--------+--------+--------+--------+--------+--------+
     |{b0, d0}|{a1, e0}|{b1, e0}|{c1, e0}|{a0, e0}|{b0, e0}|{c0, e0}|
     +--------+--------+--------+--------+--------+--------+--------+

i6 : J2G=jets(2,G); netList pack(7,edges J2G)

     +--------+--------+--------+--------+--------+--------+--------+
o7 = |{a1, c1}|{b1, c1}|{a1, d1}|{b1, d1}|{a1, e1}|{b1, e1}|{c1, e1}|
     +--------+--------+--------+--------+--------+--------+--------+
     |{c2, a0}|{d2, a0}|{e2, a0}|{c1, a0}|{d1, a0}|{e1, a0}|{c2, b0}|
     +--------+--------+--------+--------+--------+--------+--------+
     |{d2, b0}|{e2, b0}|{c1, b0}|{d1, b0}|{e1, b0}|{a2, c0}|{b2, c0}|
     +--------+--------+--------+--------+--------+--------+--------+
     |{e2, c0}|{a1, c0}|{b1, c0}|{e1, c0}|{a0, c0}|{b0, c0}|{a2, d0}|
     +--------+--------+--------+--------+--------+--------+--------+
     |{b2, d0}|{a1, d0}|{b1, d0}|{a0, d0}|{b0, d0}|{a2, e0}|{b2, e0}|
     +--------+--------+--------+--------+--------+--------+--------+
     |{c2, e0}|{a1, e0}|{b1, e0}|{c1, e0}|{a0, e0}|{b0, e0}|{c0, e0}|
     +--------+--------+--------+--------+--------+--------+--------+
\end{verbatim}

As predicted in \cite[Theorem 3.1]{2104.08933}, all jets have the same
chromatic number.

\begin{verbatim}
i8 : apply({G,J1G,J2G},chromaticNumber)

o8 = {3, 3, 3}

o8 : List
\end{verbatim}

By contrast, jets may not preserve the property of being co-chordal.

\begin{verbatim}
i9 : apply({G,J1G,J2G},x -> isChordal complementGraph x)

o9 = {true, true, false}

o9 : List
\end{verbatim}

Using Fr\"oberg's Theorem \cite{MR1171260}, we deduce that although
the edge ideal of a graph may have a linear free resolution, the edge
ideals of its jets may not have linear resolutions.

Finally, we compare minimal vertex covers of the graph and of its
second order jets.

\begin{verbatim}
i10 : vertexCovers G

o10 = {a*b*c, a*b*e, c*d*e}

o10 : List

i11 : netList pack(2,vertexCovers J2G)

      +--------------------------+--------------------------+
o11 = |a2*b2*c2*a1*b1*c1*a0*b0*c0|a2*b2*e2*a1*b1*e1*a0*b0*e0|
      +--------------------------+--------------------------+
      |a2*b2*a1*b1*c1*a0*b0*c0*e0|a2*b2*a1*b1*e1*a0*b0*c0*e0|
      +--------------------------+--------------------------+
      |c2*d2*e2*c1*d1*e1*c0*d0*e0|a1*b1*c1*a0*b0*c0*d0*e0   |
      +--------------------------+--------------------------+
      |a1*b1*e1*a0*b0*c0*d0*e0   |c1*d1*e1*a0*b0*c0*d0*e0   |
      +--------------------------+--------------------------+
\end{verbatim}

With the exception of the second row, many vertex covers arise as
indicated in \cite[Proposition 5.2, 5.3]{2104.08933}.

\section{Jets of determinantal varieties}
\label{sec:jets-determ-vari}

Determinantal varieties are classical geometric objects whose jets
have been studied to a certain degree of success
\cite{MR2100311,MR2166800,MR2389248,MR3270176,MR3020097,MR4125922}.
For our example, we consider the determinantal varieties $X_r$ of
$3\times 3$ matrices of rank at most $r$, which are defined by the
vanishing of minors of size $r+1$.  We illustrate computationally some
of the known results about jets.

\begin{verbatim}
i1 : needsPackage "Jets";

i2 : R=QQ[x_(1,1)..x_(3,3)];

i3 : G=genericMatrix(R,3,3)

o3 = | x_(1,1) x_(2,1) x_(3,1) |
     | x_(1,2) x_(2,2) x_(3,2) |
     | x_(1,3) x_(2,3) x_(3,3) |

             3       3
o3 : Matrix R  <--- R
\end{verbatim}

Since $X_0$ is a single point, its first jet scheme consists of a
single (smooth) point.

\begin{verbatim}
i4 : I1=minors(1,G);

o4 : Ideal of R

i5 : JI1=jets(1,I1);

o5 : Ideal of QQ[x0   ..x0   ][x1   ..x1   ]
                   1,1    3,3    1,1    3,3

i6 : dim JI1, isPrime JI1

o6 = (0, true)

o6 : Sequence
\end{verbatim}

The jets of $X_2$ (the variety of maximal minors) are known to be
irreducible (see \cite[Theorem 3.1]{MR2100311} or \cite[Corollary
4.13]{MR3020097}).

\begin{verbatim}
i7 : I3=minors(3,G);

o7 : Ideal of R

i8 : JI3=jets(1,I3);

o8 : Ideal of QQ[x0   ..x0   ][x1   ..x1   ]
                   1,1    3,3    1,1    3,3

i9 : isPrime JI3

o9 = true
\end{verbatim}

As for the case of $2\times 2$ minors, \cite[Theorem 5.1]{MR2100311},
\cite[Theorem 5.1]{MR2389248}, and \cite[Corollary 4.13]{MR3020097}
all count the number of components; the first two references describe
the components further. As expected, the first jet scheme of $X_1$ has
two components, one of them an affine space.

\begin{verbatim}
i10 : I2=minors(2,G);

o10 : Ideal of R

i11 : JI2=jets(1,I2);

o11 : Ideal of QQ[x0   ..x0   ][x1   ..x1   ]
                    1,1    3,3    1,1    3,3

i12 : P=primaryDecomposition JI2; #P

o13 = 2

i14 : P_1

o14 = ideal (x0   , x0   , x0   , x0   , x0   , x0   , x0   , x0   , x0   )
               3,3    3,2    3,1    2,3    2,2    2,1    1,3    1,2    1,1

o14 : Ideal of QQ[x0   ..x0   ][x1   ..x1   ]
                    1,1    3,3    1,1    3,3
\end{verbatim}

The other component is the so-called principal component of the jet
scheme, i.e., the Zariski closure of the first jets of the smooth
locus of $X_1$. To check this, we first establish that the first jet
scheme is reduced (i.e. its ideal is radical), then use the
\texttt{principalComponent} method with the option
\texttt{Saturate=>false} to speed up computations\footnote{We invite
  the reader to consult the package documentation for more details.}.

\begin{verbatim}
i15 : radical JI2==JI2

o15 = true

i16 : P_0 == principalComponent(1,I2,Saturate=>false)

o16 = true
\end{verbatim}

Finally, as observed in \cite[Theorem 18]{MR3270176}, the Hilbert
series of the principal component of the first jet scheme of $X_1$ is
the square of the Hilbert series of $X_1$.

\begin{verbatim}
i17 : apply({P_0,I2}, X -> hilbertSeries(X,Reduce=>true))

                   2     3    4            2
       1 + 8T + 18T  + 8T  + T   1 + 4T + T
o17 = {------------------------, -----------}
                      10                  5
               (1 - T)             (1 - T)

o17 : List

i18 : numerator (first oo) == (numerator last oo)^2

o18 = true
\end{verbatim}

% \bibliographystyle{halpha}
% \bibliography{\string~/Work/Bibliography/math.bib}

\def\cprime{$'$} \def\Dbar{\leavevmode\lower.6ex\hbox to 0pt{\hskip-.23ex
  \accent"16\hss}D} \def\Dbar{\leavevmode\lower.6ex\hbox to 0pt{\hskip-.23ex
  \accent"16\hss}D} \def\Dbar{\leavevmode\lower.6ex\hbox to 0pt{\hskip-.23ex
  \accent"16\hss}D} \def\Dbar{\leavevmode\lower.6ex\hbox to 0pt{\hskip-.23ex
  \accent"16\hss}D}

\end{document}